\documentclass[11pt]{article}

\usepackage{amssymb}

\newcommand{\newc}{\newcommand}

\newc{\eqnoset}{\setcounter{equation}{0}}

\newcommand{\mref}[1]{(\ref{#1})}
\newcommand{\reflemm}[1]{Lemma~\ref{#1}}
\newcommand{\refrem}[1]{Remark~\ref{#1}}
\newcommand{\reftheo}[1]{Theorem~\ref{#1}}
\newcommand{\refdef}[1]{Definition~\ref{#1}}
\newcommand{\refcoro}[1]{Corollary~\ref{#1}}
\newcommand{\refprop}[1]{Proposition~\ref{#1}}
\newcommand{\refsec}[1]{Section~\ref{#1}}

\newcommand{\beq}{\begin{equation}}
\newcommand{\eeq}{\end{equation}}
\newcommand{\beqno}[1]{\begin{equation}\label{#1}}

\newcommand{\barr}{\begin{array}}
\newcommand{\earr}{\end{array}}

\newc{\bearr}{\begin{eqnarray*}}
\newc{\eearr}{\end{eqnarray*}}

\newc{\bearrno}[1]{\begin{eqnarray}\label{#1}}
\newc{\eearrno}{\end{eqnarray}}

\newc{\non}{\nonumber}
\newc{\nol}{\nonumber\nl}

\newcommand{\bdes}{\begin{description}}
\newcommand{\edes}{\end{description}}
\newc{\benu}{\begin{enumerate}}
\newc{\eenu}{\end{enumerate}}
\newc{\btab}{\begin{tabular}}
\newc{\etab}{\end{tabular}}

\newtheorem{theorem}{Theorem}[section]
\newtheorem{defi}[theorem]{Definition}
\newtheorem{lemma}[theorem]{Lemma}
\newtheorem{rem}[theorem]{Remark}
\newtheorem{exam}[theorem]{Example}
\newtheorem{propo}[theorem]{Proposition}
\newtheorem{corol}[theorem]{Corollary}

\newcommand{\btheo}[1]{\begin{theorem}\label{#1}}
\newc{\brem}[1]{\begin{rem}\label{#1}\em}
\newc{\bexam}[1]{\begin{exam}\label{#1}\em}
\newc{\bdefi}[1]{\begin{defi}\label{#1}}
\newcommand{\blemm}[1]{\begin{lemma}\label{#1}}
\newcommand{\bprop}[1]{\begin{propo}\label{#1}}
\newcommand{\bcoro}[1]{\begin{corol}\label{#1}}
\newcommand{\etheo}{\end{theorem}}
\newcommand{\elemm}{\end{lemma}}
\newcommand{\eprop}{\end{propo}}
\newcommand{\ecoro}{\end{corol}}
\newc{\erem}{\end{rem}}
\newc{\eexam}{\end{exam}}
\newc{\edefi}{\end{defi}}

\newc{\rmk}[1]{{\bf REMARK #1: }}
\newc{\DN}[1]{{\bf DEFINITION #1: }}

\newcommand{\bproof}{{\bf Proof:~~}}
\newc{\eproof}{{\vrule height8pt width5pt depth0pt}\vspace{3mm}}

\newc{\bfrac}[2]{\dspl{\frac{#1}{#2}}}

\newc{\nid}{\noindent}

\newcommand{\dspl}{\displaystyle}
\newc{\grad}{\nabla}
\newc{\Div}{\mbox{div}}
\newc{\pdt}[1]{\dspl{\frac{\partial{#1}}{\partial t}}}
\newc{\pdn}[1]{\dspl{\frac{\partial{#1}}{\partial \nu}}}
\newc{\pdNi}[1]{\dspl{\frac{\partial{#1}}{\partial \mathcal{N}_i}}}
\newc{\pD}[2]{\dspl{\frac{\partial{#1}}{\partial #2}}}
\newc{\dt}{\dspl{\frac{d}{dt}}}
\newc{\bdry}[1]{\mbox{$\partial #1$}}
\newc{\sgn}{\mbox{sign}}

\newc{\Hess}[1]{\frac{\partial^2 #1}{\pdh z_i \pdh z_j}}
\newc{\hess}[1]{\partial^2 #1/\pdh z_i \pdh z_j}

\newc{\ag}{\alpha}
\newc{\bg}{\beta}
\newc{\cg}{\gamma}\newc{\Cg}{\Gamma}
\newc{\dg}{\delta}\newc{\Dg}{\Delta}
\newc{\eg}{\varepsilon}
\newc{\zg}{\zeta}
\newc{\thg}{\theta}
\newc{\llg}{\lambda}\newc{\LLg}{\Lambda}
\newc{\kg}{\kappa}
\newc{\rg}{\rho}
\newc{\sg}{\sigma}\newc{\Sg}{\Sigma}
\newc{\tg}{\tau}
\newc{\fg}{\phi}\newc{\Fg}{\Phi}
\newc{\vfg}{\varphi}
\newc{\og}{\omega}\newc{\Og}{\Omega}
\newc{\pdh}{\partial}

\newc{\ccG}{{\cal G}}

\newc{\ii}[1]{\int_{#1}}
\newc{\iidx}[2]{{\dspl\int_{#1}~#2~dx}}
\newc{\bii}[1]{{\dspl \ii{#1} }}
\newc{\biii}[2]{{\dspl \iii{#1}{#2} }}
\newc{\su}[2]{\sum_{#1}^{#2}}
\newc{\bsu}[2]{{\dspl \su{#1}{#2} }}

\newc{\biiom}[1]{{\dspl\int_{\bdrom}~ #1 ~d\sg}}
\newc{\io}[1]{{\dspl\int_{\Og}~ #1 ~dx}}
\newc{\bio}[1]{{\dspl\int_{\bdrom}~ #1 ~d\sg}}
\newc{\bsir}{\bsu{i=1}{r}}
\newc{\bsim}{\bsu{i=1}{m}}

\newc{\iibr}[2]{\iidx{\bprw{#1}}{#2}}
\newc{\Intbr}[1]{\iibr{R}{#1}}
\newc{\intbr}[1]{\iibr{\rg}{#1}}
\newc{\intt}[3]{\int_{#1}^{#2}\int_\Og~#3~dxdt}

\newc{\itQ}[2]{\dspl{\int\hspace{-2.5mm}\int_{#1}~#2~dz}}
\newc{\mitQ}[2]{\dspl{\rule[1mm]{4mm}{.3mm}\hspace{-5.3mm}\int\hspace{-2.5mm}\int_{#1}~#2~dz}}
\newc{\mitQQ}[3]{\dspl{\rule[1mm]{4mm}{.3mm}\hspace{-5.3mm}\int\hspace{-2.5mm}\int_{#1}~#2~#3}}

\newc{\mitx}[2]{\dspl{\rule[1mm]{3mm}{.3mm}\hspace{-4mm}\int_{#1}~#2~dx}}
\newc{\mitmu}[2]{\dspl{\rule[1mm]{3mm}{.3mm}\hspace{-4mm}\int_{#1}~#2~d\mu}}
\newc{\iidmu}[2]{{\dspl\int_{#1}~#2~d\mu}}

\newc{\iidm}[3]{{\dspl\int_{#1}~#2~d #3}}

\newc{\itQmu}[2]{\dspl{\int\hspace{-2.5mm}\int_{#1}~#2~d\mu}}
\newc{\mitQmu}[2]{\dspl{\rule[1mm]{4mm}{.3mm}\hspace{-5.3mm}\int\hspace{-2.5mm}\int_{#1}~#2~d\mu}}

\newc{\mitQq}[2]{\dspl{\rule[1mm]{4mm}{.3mm}\hspace{-5.3mm}\int\hspace{-2.5mm}\int_{#1}~#2~d\bar{z}}}
\newc{\itQq}[2]{\dspl{\int\hspace{-2.5mm}\int_{#1}~#2~d\bar{z}}}

\newc{\pder}[2]{\dspl{\frac{\partial #1}{\partial #2}}}

\newc{\bdrom}{\bdry{\Og}}

\newc{\bilhom}{\mbox{Bil}(\mbox{Hom}(\RR^{nm},\RR^{nm}))}
\newc{\VV}[1]{{V(Q_{#1})}}

\newc{\ccA}{{\mathcal A}}
\newc{\ccB}{{\mathcal B}}
\newc{\ccC}{{\mathcal C}}
\newc{\ccD}{{\mathcal D}}
\newc{\ccE}{{\mathcal E}}
\newc{\ccH}{\mathcal{H}}
\newc{\ccF}{\mathcal{F}}
\newc{\ccI}{{\mathcal I}}
\newc{\ccJ}{{\mathcal J}}
\newc{\ccK}{{\mathcal K}}
\newc{\ccP}{{\mathcal P}}
\newc{\ccQ}{{\mathcal Q}}
\newc{\ccR}{{\mathcal R}}
\newc{\ccS}{{\mathcal S}}
\newc{\ccT}{{\mathcal T}}
\newc{\ccX}{{\mathcal X}}
\newc{\ccY}{{\mathcal Y}}
\newc{\ccZ}{{\mathcal Z}}

\newc{\bb}[1]{{\mathbf #1}}

\newc{\myprod}[1]{\langle #1 \rangle}
\newc{\mypar}[1]{\left( #1 \right)}

\newc{\BLLg}{\mathbf{\LLg}}

\newc{\mA}{\mathbf{A}}
\newc{\mB}{\mathbf{B}}
\newc{\mC}{\mathbf{C}}
\newc{\mD}{\mathbf{D}}
\newc{\mE}{\mathbf{E}}
\newc{\mF}{\mathbf{F}}
\newc{\mJ}{\mathbf{J}}
\newc{\mG}{\mathbf{G}}
\newc{\mP}{\mathbf{P}}
\newc{\mR}{\mathbf{R}}
\newc{\mQ}{\mathbf{Q}}
\newc{\mX}{\mathbf{X}}
\newc{\muu}{\mathbf{u}}
\newc{\mvv}{\mathbf{v}}

\newc{\mllg}{\mathbb{\lambda}}
\newc{\mLLg}{\mathbf{\LLg}}

\newc{\lspn}[2]{\mbox{$\| #1\|_{\Lsp{#2}}$}}
\newc{\Lpn}[2]{\mbox{$\| #1\|_{#2}$}}
\newc{\Hn}[1]{\mbox{$\| #1\|_{H^1(\Og)}$}}

\newc{\mynorm}[2]{\| #1\|_{#2}}

\newcommand{\RR}{{\rm I\kern -1.6pt{\rm R}}}

\newc{\itQQ}[2]{\dspl{\int_{#1}#2\,dz}}
\newc{\mmitQQ}[2]{\dspl{\rule[1mm]{4mm}{.3mm}\hspace{-4.3mm}\int_{#1}~#2~dz}}
\newc{\MmitQQ}[2]{\dspl{\rule[1mm]{4mm}{.3mm}\hspace{-4.3mm}\int_{#1}~#2~d\mu}}

\newc{\MUmitQQ}[3]{\dspl{\rule[1mm]{4mm}{.3mm}\hspace{-4.3mm}\int_{#1}~#2~d#3}}
\newc{\MUitQQ}[3]{\dspl{\int_{#1}~#2~d#3}}

\newc{\mccP}{\mathbb{P}}
\newc{\mccK}{\mathbb{K}}

\newc{\DKTmU}{\mccK(U)}
\newc{\DKTmUold}{(K_U(U)^{-1})^T}

\newc{\myPi}{\mathbf{W}}
\newc{\myIbar}{\bar{\ccI}_1}
\newc{\myIhat}{\hat{\ccI}_1}
\newc{\myIbreve}{\breve{\ccI}_0}

\newc{\mmk}{\mathbf{k}}

\newcommand{\ma}{\mathbf{a}}
\newcommand{\mg}{\mathbf{g}}

\newc{\mfu}{\mathbf{f_u}}
\newc{\mh}{\mathbf{h}}
\newcommand{\barrl}[2]{\barr{ll}\lefteqn{#1}\hspace{#2}&\\}

\begin{document}

\vspace*{-.8in}
\begin{center} {\LARGE\em Uniqueness and Regularity of Unbounded Weak Solutions to a Class of Cross Diffusion Systems.}

 \end{center}

\vspace{.1in}

\begin{center}

{\sc Dung Le}{\footnote {Department of Mathematics, University of
Texas at San
Antonio, One UTSA Circle, San Antonio, TX 78249. {\tt Email: Dung.Le@utsa.edu}\\
{\em
Mathematics Subject Classifications:} 35J70, 35B65, 42B37.
\hfil\break\indent {\em Key words:} Cross diffusion systems,  H\"older
regularity, global existence.}}

\end{center}

\begin{abstract}
We establish the uniqueness and regularity of weak (and very weak) solutions to a class of cross diffusion systems which is  inspired by models in mathematical biology/ecology, in particular the Shigesada-Kawasaki-Teramoto (SKT) model in population biology. No boundedness assumption on these solutions is supposed here as known techniques for scalar equations such as maximum/comparison principles are generally unavailable for systems.  Furthermore, for planar domains  we show that unbounded weak solutions satisfying mild integrability conditions are in fact smooth.   \end{abstract}

\vspace{.2in}

\section{Introduction}\label{introsec}\eqnoset

In this paper, we study the following parabolic system of $m$ equations ($m\ge2$)  for the unknown vector $u=[u_i]_{i=1}^m$
\beqno{ep1}u_t=\Delta(P(u))+f(u),\quad (x,t)\in Q:=\Og\times(0,T_0).\eeq  
Here,  $P:\RR^m\to\RR^m$ is a $C^2$ map and $f:\RR^m\to \RR^m$ is a $C^1$ map. $\Og$ is a bounded domain with smooth boundary in $\RR^N$, $N\ge2$, and $T_0>0$. 

The system is equipped with boundary and initial conditions
\beqno{ep1bc}\left\{\barr{l} \mbox{$u=0$  on $\partial \Og\times(0,T_0)$},\\ u(x,0)=u_0(x) ,\quad x\in \Og. \earr\right.\eeq

The consideration of \mref{ep1} is motivated by the extensively studied  model in population biology introduced by Shigesada {\it et al.} in \cite{SKT}
\beqno{e0}\left\{\barr{lll} u_t &=& \Delta(d_1u+\ag_{11}u^2+\ag_{12}uv)+k_1u+\bg_{11}u^2+\bg_{12}uv,\\v_t &=& \Delta(d_2v+\ag_{21}uv+\ag_{22}v^2)+k_2v+\bg_{21}uv+\bg_{22}v^2.\earr\right.\eeq   Here, $d_i,\ag_{ij},\bg_{ij}$ and $k_i$ are constants with $d_i>0$. Dirichlet or Neumann boundary conditions were usually assumed for \mref{e0}. This model was used to describe the population dynamics of {\em two} species densities  $u,v$ which move and react under the influence of population pressures.

Under suitable assumptions on the constant parameters $\ag_{ij}$'s, $\bg_{ij}$'s and that $\Og$ is a planar domain ($N=2$), Yagi proved in \cite{yag} the global existence of (strong) positive solutions, with positive initial data. In this paper, among other general settings, we will investigate weak solutions to multi-species versions of \mref{e0} for more than two species and consider much more general  structural conditions. Naturally, we will replace the quadratics in the Laplacians and $f_i$ of \mref{e0} by polynomials of order $k+1$ for some $k>0$. Obviously,
the system \mref{e0} is a special case of \mref{ep1} with $m=2$ and 
$P:\RR^2\to\RR^2$ being a quadratic map.

Let us describe the multi-species version of \mref{e0}.
For $m\ge2$ let $u:\Og\times(0,T_0)\to \RR^m$ and  $\ag_i,\bg_i\in \RR^m$, $i=1,\ldots,m$. The multi-species version of \mref{e0} is then \mref{ep1} with  \beqno{sktpf}P_i(u)=d_iu_i+u_i\myprod{\ag_i,u},\; f_i(u)=k_iu_i+u_i\myprod{\bg_i,u}.\eeq

The generalized multi-species version of \mref{e0} is naturally obtained by replacing $\ag_i, \bg_i$ by maps on $\RR^m$ with certain polynomial growth. That is, $|\ag_i(u)|,|\bg_i(u)|\sim|u|^\kappa$ for $u\in \RR^m$ and some $\kappa\ge 0$. 
In this paper, with this setting, we will refer to \mref{ep1} as the SKT system if $\kappa=0$ and the generalized SKT one if $\kappa>0$. We should mention a recent work \cite{LM} where a similar setting was considered.

We can write \mref{ep1} in its divergence form $$ u_t=\Div(A(u)Du) + f(u),$$ where $A(u)=P_u(u)$, the Jacobian of $P(u)$. For appropriate $\ag_{ij}$'s (see \cite{yag} for \mref{e0}) this system is parabolic in the sense that $P_u$ is elliptic. That is, there are  a constant $\llg_0>0$ and a {\em linear} function $\llg(u)= \llg_0+|u|$ such that $\llg (u)\ge\llg_0$ and $$ \myprod{P_u(u)\zeta,\zeta}\ge \llg(u) |\zeta|^2, \mbox{ for all $u\in\RR^m, \zeta\in\RR^{Nm}$}.$$ This will be the structural main assumption on \mref{ep1} and we note that, in this paper, the ellipticity function $\llg(u)$ needs not to be of polynomial growth in general.

Our first goal is to establish the uniqueness of {\em unbounded} weak solutions of \mref{ep1}. To our best knowledge, this is the first time this problem is treated for cross diffusion {\em systems} like \mref{ep1} with no boundedness is assumed. The problem was addressed in \cite{BC} for scalar equations, where maximum principles were available, and it was nontrivial already. Here, we are working with cross diffusion systems and no invariant/maximum principles are available so that the uniqueness question must be treated in a completely different way. Also, without boundedness, the matrix $A(u)=P_u(u)$ is not regular elliptic as it is not bounded from above.
Under very weak integrability conditions on $u$, we obtain the uniqueness of weak (and very weak) solution of \mref{ep1} in \reftheo{uniweak}. Its immediate consequences then apply to the generalized weak solutions from $V_2(Q)$ (following the common definition of \cite{LSU} in literature) of the SKT system \mref{e0} and its generalized versions.

Next, we consider the regularity properties of weak solutions (again, unbounded). This is a long standing and hard problem in the theory of pde's, especially for strongly coupled parabolic systems like \mref{ep1} or even \mref{e0}. Here, by combining our uniqueness results with the theory in \cite{dleANS,dlebook} (see also \cite{letrans}) which dealt with strong solutions, we consider \mref{ep1} defined on planar domains ($N=2$) and show that {\em unbounded} generalized solutions from $V_2(Q)$ of the SKT systems and its generalized versions are in fact classical. This result is new and our indirect approach (a combination of the studies of uniqueness of weak solutions and existence of strong ones) may come as a surprise in comparison with direct methods which worked only with {\em bounded} weak solutions but not in our case (see \refrem{regrem}).

The paper is organized as follows. In \refsec{mainres}, we state the general structural conditions on \mref{ep1}, the main uniqueness and regularity results  and some of their immediate applications to the SKT systems. More examples will be discussed later in \refsec{uniboundw} and \refsec{regsec} when we complete the proof of the main results because some improvements can be obtained by slight modifications of the proof thanks to additional and special structures of the models and need more detailed discussion. Technical tools will be presented in \refsec{compsec}. \refsec{uniboundw} is devoted to the proof of the uniqueness result. We  prove and discuss the regularity results in \refsec{regsec}.

\section{The main results}\eqnoset\label{mainres}
We state our main results in this section. We first discuss the uniqueness of weak solutions. This has been done for scalar parabolic equations for {\em bounded} weak solutions. But this is not the case for systems because the boundedness of solutions to systems generally is an open problem and the arguments for scalar equations are not applicable here. For the system \mref{ep1} we will establish this result for {\em unbounded} weak solutions which is defined in a very general sense and satisfy very mild integrability conditions.

The system \mref{ep1} with boundary and initial condition \mref{ep1bc} is a special case of the parabolic system in divergent form with $A(u)=P_u(u)$, the Jacobian of $P(u)$,
\beqno{ep1div} u_t=\Div(A(u)Du) + f(u).\eeq
Following the standard definition, we say that 
\bdefi{wsolndef} $u$ is a weak solution on $\Og \times(0,T_0)$ the system \mref{ep1div} with boundary and initial condition \mref{ep1bc} if $u\in L^\infty((0,T_0), L^1(\Og))$ and $A(u)Du\in L^1(\Og \times(0,T_0))$ and
for  a.e. $T\in(0,T_0)$ and any $\fg\in C^1(\Og \times (0,T))$ we have
\beqno{wdef}\barrl{\iidx{\Og}{\myprod{u(T),\fg(T)}-\myprod{u_0,\fg(0)}}=}{3cm}&\itQ{\Og \times(0,T)}{[\myprod{u, \fg_t}-\myprod{A(u)Du,D\fg} +\myprod{f(u),\fg}]}.\earr\eeq
\edefi

In this paper, as $A(u)Du=D(P(u))$, from  the equation \mref{wdef} with a simple integration by parts in $x$, assuming that $P(u)=0$ on $\partial\Og\times(0,T_0)$ and using the homogeneous Dirichlet boundary condition, we also have the following weaker definition for weak solution $u$ of \mref{ep1} with no integrability assumption on (weak) derivatives of $u$ and more restrictive admissible test functions.
\bdefi{wsolndef0} $u$ is a {\em very} weak solution on $\Og \times(0,T_0)$ the system \mref{ep1} with boundary and initial condition \mref{ep1bc} if $u\in L^\infty((0,T_0), L^1(\Og))$ and $P(u)\in L^1(\Og \times(0,T_0))$ and
for  a.e. $T\in(0,T_0)$ and any $\fg\in C^2(\Og \times (0,T))$ we have
\beqno{Pwdef0}\barrl{\iidx{\Og}{\myprod{u(T),\fg(T)}-\myprod{u_0,\fg(0)}}=}{3cm}&\itQ{\Og \times(0,T)}{[\myprod{u, \fg_t}+\myprod{P(u),\Delta\fg} + \myprod{f(u),\fg}]}.\earr\eeq
\edefi

We note that a admissible test function $\fg$ in this definition, as a minimum requirement, needs only that $\fg_t, \Delta \fg\in L^\infty(\Og \times(0,T_0))$. Also, in both definitions, we can consider initial data $u_0\in L^1(\Og)$.

In order for \mref{ep1div} to be regular parabolic with $A(u)=P_u(u)$, we naturally impose our main assumption on the structure of the system is
\bdes\item[A)] $P(0)=0$. $P_u$ is regular elliptic. That is there are function $\llg$ and constant $\llg_0>0$ such that $\llg (u)\ge\llg_0$ and $$ \myprod{P_u(u)\zeta,\zeta}\ge \llg(u) |\zeta|^2, \mbox{ for all $u\in\RR^m, \zeta\in\RR^{Nm}$}.$$
\edes
Furthermore, 
\bdes\item[F)]there is a {\em convex} function $\hat{F}$ such that $|\partial_u f(u)|^2\llg (u)^{-1}\le \hat{F}(u)$ on $\RR^m$.\edes

We also introduce the following notations for our theorem statements.
\beqno{pr} p_r:=\frac{r'p}{p-r'}, \mbox{ where $r'=r/(1-r)$, the conjugate of $r$},\eeq
\beqno{qPsidef}\sg_N=\left\{\barr{ll} \mbox{any number in $(1,\infty)$} & \mbox{ if $N=2 $},\\ \mbox{any number in $(1,6+\frac{10}{3})$} & \mbox{ if $N= 3$},\\ \frac{2(N+2)}{N-2} & \mbox{ if $N\ge 4$}.\earr\right.\eeq

The main result of this paper is the following uniqueness theorem.

\btheo{uniweak} Assume A) and F). For some $p>2$ we also assume the following integrability continuity conditions (with $Q=\Og\times(0,T_0)$ and the notations in \mref{pr} and \mref{qPsidef}):
\bdes \item[i)]
The map $u\to \partial_u P(u)$ is continuous from $L^p(Q)$ to $L^{p_2}(Q)$,
\item[ii)] The map  $u\to \partial_u f(u)$ is continuous from $L^p(Q)$ to $L^{p_{\sg_N}}(Q)$.
\edes
If $u$ is a {\em very} weak solution, in the sense of \refdef{wsolndef0}, and satisfies $u\in L^p(Q)$ and for some $q_0\ge N/2$ \beqno{llguw0} \sup_{t\in (0,T_0)}\|\hat{F}(u(t))\|_{L^{q_0}(\Og)}<\infty,\eeq then $u$ is unique. \etheo

We will see that the conditions of this theorem can be verified in many models in application. To discuss this matter further, we note that there are many ways to define the concept of weak solution and they all start with the equation \mref{wdef} in \refdef{wsolndef} (or \mref{Pwdef0} in \refdef{wsolndef0}) and there is a trade off among these definitions in order that the integrals of \mref{wdef} are all finite. The main difference lies in the choices  of  admissible test function $\fg$ in \mref{wdef}.  If the space of admissible test functions is more restrictive then the space of weak solutions will be in wider and it is harder to obtain the uniqueness result and vice versa. We would like to remark this fact here for future references.

\brem{wPdefrem} Of couse, \refdef{wsolndef0} is weaker than \refdef{wsolndef} so that \reftheo{uniweak} also applies to weak solutions in the sense of \refdef{wsolndef}. \erem

Still, our \refdef{wsolndef} is an enough general one as we needs only that the first order derivatives of $\fg $ are defined and $\fg \in C^1(Q)$. Consequently, a weak solution $u$ is this sense needs only satisfy $u\in L^\infty((0,T_0),L^1(\Og))$ and $D(P(u))$, $f(u)$ are in $L^1(Q)$ in order that the integrals in \mref{wdef} are all finite. Of course, our class of weak solutions is sufficiently wide and the checking of their integrability conditions of \reftheo{uniweak} seems to be already not an easy condition under such limited information.

On the other hand, if we allow more general test function $\fg$ then the space of weak solutions will be smaller and the uniqueness result can be applied easily and almost immediate in some cases. 

Following \cite[Chapter III]{LSU}, which has been used widely in literature, we say that $u$ is a generalized solution from $V_2(Q)$, the  Banach space with norm
$$\|u\|_{V_2(Q)}=\sup_{t\in(0,T_0)}\|u\|_{L^2(\Og\times \{t\})}+\|Du\|_{L^2(Q)},$$
if $u$ satisfies \mref{wdef} for any test function $\fg \in W_2^{1,1}(Q)$, the Hilbert space with scalar product $$ \myprod{u,v}_{W_2^{1,1}(Q)}=\itQ{Q}{[\myprod{u,v}+\myprod{u_t,v_t}+\myprod{Du,Dv}]}.$$

Adopting this concept of weak solutions, we have
\bcoro{SKTuni} Generalized solutions from $V_2(Q)$ to the SKT system  on domains in $\RR^N$ with $N\le 4$ are unique. \ecoro

Much more  on the generalized version of the SKT systems will be discussed in \refsec{uniboundw}.

Next, we will consider the regularity of {\em unbounded} weak (and very weak) solutions of \mref{ep1}. This is a very hard problem in the theory of pdes, especially for strongly coupled parabolic systems with \mref{ep1} as a special case. There is a vast literature on this problem, see \cite{GiaS}, and all assume that the considered weak solution is {\em bounded} and satisfies a crucial condition that its BMO norm is small in small balls. Here, for any ball $B_R$ of radius $R$ in $\RR^N$ and $\Og_R=B_R\cap\Og$ the BMO norm of a (vector valued) function $u$ is defined by
\beqno{bmodef} \|u\|_{BMO(\Og_R)}:=\sup_{B_r\subset\Og_R}\mitx{B_r}{|u-u_r|}+\iidx{\Og_R}{|u|},\eeq where $u_r$ is, as usual, the average of $u$ over $B_r$. Beside the boundedness assumption, which is already a very hard problem for weak solutions to cross diffusion systems, the smallness of the norm \mref{bmodef} when $R$ is small is even a harder problem and, to the best of our knowledge, none were done in literature to address either question in this general setting.

In this paper, we follow an indirect and novel approach to this regularity problem and the idea is very simple: If we can show that there exists a strong solution $u$ to \mref{ep1} then this solution is of course also a weak one and satisfies the integrability conditions of \reftheo{uniweak}. By the uniqueness result for such weak solutions, any weak solution satisfying sufficient integrability of \reftheo{uniweak} is exactly this strong solution $u$ and therefore it is in fact classical (or smooth).
Thus, we immediately have the following statement.
\bcoro{wsequal2} Assume that \mref{ep1} possesses a strong solution. If $u$ is a ({\em very}) weak solution of \mref{ep1}  and satisfies the integrability conditions of \reftheo{uniweak} then $u$ is a classical one.
\ecoro

Of course, the existence of a strong solution to \mref{ep1} is also one of the hardest problems in the theory of cross diffusion systems so that the hypothesis of \refcoro{wsequal2} is a very bold one. This existence problem was considered in the pioneering work by Amann \cite{Am2} (see also \cite{yag}), using semigroup theory. However, his hypothesis on  a priori estimates for the gradients of (strong) solutions was also equivalent to the investigation of H\"older continuity of the solutions, again a difficult task for systems and only few works were done and relied on very ad hoc techniques which applied only to special cases of \mref{ep1}. Recently, in \cite{dlebook}, we introduced an alternative approach using fixed point theories to establish the existence of strong solutions under a set of a priori integrability conditions and, again crucially, the smallness of the norm \mref{bmodef}.
One of the advantages of the theory is that one can remove the boundedness assumption in
\cite{Am2} and replace it by certain mild integrability ones. 

Yet, establishing the smallness of the norm \mref{bmodef} eludes many efforts for general dimension $N$.
We present in \refsec{regsec} the conditions S) and S') which can be verified for \mref{ep1} in applications to affirm both integrability and smallness of the norm \mref{bmodef} when $N=2$.
We restrain ourselves from giving the details of these conditions as they need some technical preparations and discussion. Here, we just state the following consequence for the SKT systems and its generalizations. 
\bcoro{SKTrek}  Suppose further that the reaction term $f(u)$ satisfies \beqno{SKTfu} \myprod{f(w),w}\le \eg_0\llg(w)|w|^2+ C|w|^2 \eeq for some positive constants $C,\eg_0$. If $\eg_0$ is sufficiently small, in terms of the diameter of $\Og$, then the generalized solution from $V_2(Q)$ of  of the SKT system and its generalized versions (for $k<2$) on planar domains are classical.
\ecoro

Even in this special case ($N=2$), the result is new and remarkable because no boundedness is assumed.

\section{Some technical lemmas} \eqnoset \label{compsec}

In this section, we collect some technical lemmas some of which may be elementary to experts and others are subtle  will play crucial roles in the proof of our main results and examples.
 
First of all, in the proof we will frequently make use of the following interpolation Sobolev inequality 
 
 \blemm{Sobointineq} For any  $\eg>0$, $\bg\in(0,1]$, $p\ge1$ and $W\in W^{1,p}(\Og)$ we can find a constant $C(\eg,\bg)$  such that
 \beqno{intineq}\|W\|_{L^q(\Og)}\le \eg\|DW\|_{L^p(\Og)} + C(\eg,\bg)\|W^\bg\|_{L^1(\Og)}^\frac{1}{\bg} \mbox{ for any $q\in[1,p_*)$}.\eeq \elemm
 
 Here and throughout this paper, $p_*$ as usual, denotes the Sobolev conjugate of $p$.
 $$p_*=\left\{\barr{ll} Np/(N-p)& \mbox{if $p<N$},\\
 \mbox{any number in $(1,\infty)$}& \mbox{otherwise}.\earr \right.$$
 
\bproof By contradiction, assume that \mref{intineq} is not true then we can find $\eg_0>0$ and a sequence $\{W_n\}$ such that \beqno{intineqX}\|W_n\|_{L^q(\Og)}> \eg_0\|DW_n\|_{L^p(\Og)} + n\|W_n^\bg\|_{L^1(\Og)}^\frac{1}{\bg} \mbox{ for any $n$}.\eeq By scaling we can suppose that $\|W_n\|_{L^q(\Og)}=1$. The above implies that $\|DW_n\|_{L^p(\Og)}<1/\eg_0$ for all $n$. We see that $\{W_n\}$ is bounded in $W^{1,p}(\Og)$ so that, by compactness (as $q<p_*$), we can assume that it converges to some $W$ in $L^q(\Og)$. Of course, $\|W\|_{L^q(\Og)}=1$. Meanwhile, \mref{intineqX} implies $\|W_n^\bg\|_{L^1(\Og)}\to0$ so that $\|W^\bg\|_{L^1(\Og)}=0$, which can be easily seen by H\"older's inequality and the H\"older continuity of the function $|x|^\bg $. Thus $W=0$ a.e. on $\Og $ and contradicts the fact that $\|W\|_{L^q(\Og)}=1$. The proof is complete. \eproof

Next, we need some estimates of solutions to the following linear parabolic system 
\beqno{Psidefa0}\left\{\barr{ll}\Psi_t=\ccA\Delta \Psi+\ccG \Psi&\mbox{on $Q=\Og \times(0,T)$},\\\Psi=0&\mbox{on $\partial\Og \times(0,T)$},\\ \Psi(x,0)=\psi(x).&\earr\right.\eeq

Here, $\ccA(x,t) , \ccG(x,t)$ are matrices of sizes $m\times m$ and $m\times 1$ respectively and satisfy the following condition

\bdes \item[AG)] $\ccA ,\ccG$   are smooth in $Q$ and $\ccA$ is regular elliptic. That is, there are function $\llg_*$ on $Q$ and a constant $\llg_0>0$ such that $\llg_{*}(x,t)\ge \llg_0$ and
\beqno{ccAcond} \llg_*(x,t)|\zeta|^2\le \myprod{\ccA (x,t)\zeta,\zeta} \mbox{ for all $\zeta\in \RR^{m}$ and all $(x,t)\in Q$.}\eeq 

Furthermore, $\|\ccA(\cdot,t) \|_{L^\infty(\Og)}$ is continuous at $t=0$ and the $L^\infty(\Og)$ norms of $D\ccA(\cdot,t)$ and  $\ccG(\cdot,t)$ are bounded near $t=0$. 
\edes

Next, we present a lemma which concerns the behavior of $\Psi$ near $t=0$ and is in the spirit of the Hille-Yoshida theorem (e.g., see \cite[Theorem 7.8]{brezis}) on the continuity of the (weak) derivatives of weak solutions to \mref{Psidefa0} when $\ccA $ is a constant. Otherwise, this theorem does not apply directly to our case because $\ccA$ depends on both $x,t$. However, later on we just need an appropriate bound for $\liminf_{t\to0}\|D\Psi\|_{L^2(\Og \times\{t\})}$ for our purposes in this paper and one should also note that this bound does not depend on $\|\ccG_u\|_{L^\infty(Q)}$ and $\|D\ccA_u\|_{L^\infty(Q)}$. As we cannot find an appropriate reference, we state the fact and present its proof in details.

\blemm{limDuat0} Assume AG) and $\psi\in W^{1,2}(\Og)$. 
If $\Psi$ is a strong solution  to \mref{Psidefa0} then \beqno{limDPsi}\liminf_{t\to0}\|D\Psi\|_{L^2(\Og \times\{t\})}\le \|D\psi\|_{L^2(\Og)}.\eeq
\elemm

\bproof We split $\Psi=h+H$ with $h,H$ solving 
$$h_t=a(x) \Delta h +\ccG h,\; H_t= \ccA \Delta H +\ccB \Delta h+ \ccG H,$$ where  $a(x)=\ccA(x,0)$ and $\ccB(x,t)=\ccA(x,t)-a(x)$. Also, $h(0)=\psi$ and $H(0)=0$. We rewrite the equation for $H$ as
$$H_t= \Div(\ccA DH +\ccB Dh)- D\ccA DH-D\ccB Dh + \ccG H$$
and test the system with $H$ and use the fact that $H(0)=0$ to obtain  for any $s>0$ that \beqno{HHH}\barrl{\iidx{\Og \times\{s\}}{|H|^2}+\int_0^s\iidx{\Og}{\myprod{\ccA DH,DH}}=}{2cm}&-\dspl{\int_0^s}\iidx{\Og}{(\myprod{\ccB Dh,DH}+\myprod{D\ccA DH+D\ccB Dh,H}+\myprod{\ccG H,H})}.\earr\eeq  Applying Young's inequalities to the integrands involving $DH$ on the right hand side  and using  the ellipticity assumption \mref{ccAcond}, we easily get
$$\barrl{\int_0^s\iidx{\Og}{|DH|^2}dt \le C\int_0^s\iidx{\Og}{|\ccB|^2 |Dh|^{2}}dt+}{3cm}& C\dspl{\int_0^s}\iidx{\Og}{(|D\ccA|^2+|\ccG|)|H|^2}dt+C\dspl{\int_0^s}\iidx{\Og}{|D\ccB|Dh|H|}dt.\earr$$
We now divide the about inequality by $s$ to have
\beqno{keyDuat0}\barrl{\frac1s\int_0^s\iidx{\Og}{|DH|^2}dt \le C\frac1s\int_0^s\iidx{\Og}{|\ccB|^2|Dh|^{2}}dt+}{2cm}& C \dspl{\frac1s\int_0^s}\iidx{\Og}{(|D\ccA|^2+|\ccG|) |H|^2}dt+C\frac1s\dspl{\int_0^s}\iidx{\Og}{|D\ccB|Dh|H|}dt.\earr\eeq
We will let $s\to0$ and need to investigate the limits of the terms on the right hand side.

From the definition of $\ccB$, $\ccB(x,t)=\ccA(x,t)-\ccA(x,0)$. As $ \|\ccA(\cdot,t) \|_{L^\infty(\Og)}$ is continuous at $t=0$, we see that $\lim_{t\to0}\|\ccB(\cdot,t)\|_{L^\infty(\Og)}=0$. Also, by the Hille-Yoshida theorem (e.g., see \cite[Theorem 7.8]{brezis}, note that $a(x)=\ccA (x,0)$ is elliptic, smooth and independent of $t$), $Dh$ belongs to $C([0,T_0],L^{2}(\Og))$. In particular, $\|Dh(t)\|_{L^{2}(\Og)}$ is  continuous at $t=0$. Hence, the limit of the first term on the right hand side of \mref{keyDuat0} when $s\to0$ is zero.

Meanwhile,  as $H$ is a strong solution and   $H(0)=0$, we have $$\frac1s\int_0^s\iidx{\Og}{|H|^2}dt\to 0.$$ By the assumption AG),  the $L^\infty(\Og)$ norms of $D\ccA(\cdot,t)$ (and then $D\ccB(\cdot,t)$) and  $\ccG(\cdot,t)$ are bounded near $t=0$. We see that the last two terms tend to $0$. 
Hence, letting $s\to0$ in \mref{keyDuat0}, we derive $\liminf_{t\to0}\|DH\|_{L^2(\Og \times\{t\})}=0$. Meanwhile, as we explained above, $\lim_{t\to0}\|Dh(t)\|_{L^{2}(\Og)}=\|D\psi\|_{L^{2}(\Og)}$ Because $\Psi=h+H$, we obtain \mref{limDPsi} of the lemma. 
\eproof

\brem{LpDADBrem} \reflemm{limDuat0} still holds if the $L^p(\Og)$ norms of $D\ccA(\cdot,t)$  and  $\ccG(\cdot,t)$ are bounded near $t=0$ for some finite $p>1$. Indeed, as $H$ is a strong solution and   $H(0)=0$, we have $$\frac1s\int_0^s\iidx{\Og}{|H|^q}dt\to 0$$  for any $q$ as $s\to0$. It is easy to see that a simple use of H\"older's inequality shows that the last two integrals of \mref{keyDuat0} still tend to $0$. 
\erem

\brem{LpDADBrem0}
The assertion of \reflemm{limDuat0} also applies to the system $\Psi_t =\Div(\ccA D\Psi)+\ccG\Psi$ which is in divergence form. The systems for $h,H$ now are
$$h_t=\Div(a(x) D h) +\ccG h,\; H_t= \Div(\ccA DH) +\Div(\ccB Dh)+ \ccG H.$$
We then repeat the same argument and obtain \mref{keyDuat0} without the terms involving $D\ccA, D\ccB$ so that the assumption on the boundedness of $D\ccA $ (and so $D\ccB$) near $t=0$ {\em can be dropped}.

\erem

We then have the following estimates for derivatives of $\Psi$ which will play an important role in our proof of uniqueness results. 

\blemm{Psiexistence} Assume AG) and let $\psi\in C^{1}(\Og,\RR^m)$.
Then there is a classical solution $\Psi$ to
\mref{Psidefa0}.
In addition, assume that there are constants $C_0,q_0$ such that $q_0\ge N/2$ ($q_0>1$ if $N=2$) and   \beqno{g*bound}\sup_{(0,T)}\|g_*\|_{L^{q_0}(\Og\times\{\tau\})}\le C_0,\; \mbox{where $g_*:=|\ccG|^2\llg_{*}^{-1}$}.\eeq 

Then there is a constant $C(T,C_0)$  such that
\beqno{Psibounda} \iidx{\Og\times\{s\}}{|D\Psi|^2} \le C(T,C_0)\|D\psi\|_{L^2(\Og)} \mbox{  for all $s\in[0,T]$},\eeq
\beqno{D2Psibound} \itQ{Q}{|\Delta\Psi|^2} \le \llg_0^{-1}C(T,C_0)\|D\psi\|_{L^2(\Og)}.\eeq
\elemm

One should note that the constant $C(T,T_0)$ in \mref{Psibounda} and \mref{D2Psibound} is independent of the norms of $\ccA ,\ccG $ and their derivatives.

\bproof The existence of $\Psi$ is clear because $\ccA , \ccG $ are smooth and $\ccA$ is regular elliptic. We just need to establish \mref{Psibounda} and \mref{D2Psibound}. Multiplying \mref{Psidefa0} with $\Delta\Psi$ and integrating by parts ($\Psi_t=0$ on the boundary as homogeneous Dirichlet condition is considered for $\Psi$), we get for any $0<s<T$ and $Q^{(s)}=\Og\times(0,s)$ $$\itQ{Q^{(s)}}{\frac{d}{dt}|D\Psi|^2}+\itQ{Q^{(s)}}{\myprod{\ccA\Delta\Psi,\Delta\Psi}}=-\itQ{Q^{(s)}}{\myprod{\ccG\Psi,\Delta\Psi}}.$$

By the ellipticity of $\ccA $ we get for any $s'\in(0,s)$
\beqno{Psistart} \iidx{\Og\times\{s\}}{|D\Psi|^2}+\itQ{Q^{(s)}}{\llg_*|\Delta\Psi|^2}\le \iidx{\Og\times\{s'\}}{|D\Psi|^2}+\itQ{Q^{(s)}}{\myprod{\ccG\Psi,\Delta\Psi}}.\eeq

Next,  by Young's inequality   
$$\myprod{\ccG\Psi,\Delta\Psi}\le \eg\llg_*|\Delta\Psi|^2+C(\eg)g_*|\Psi|^2, \mbox{ $g_*:=|\ccG|^2\llg_{*}^{-1}$}.$$ 

Therefore, for small $\eg>0$ we deduce from \mref{Psistart} the following inequality 
\beqno{Psistartz}\iidx{\Og\times\{s\}}{|D\Psi|^2}+\itQ{Q^{(s)}}{\llg_*|\Delta\Psi|^2}\le \iidx{\Og\times\{s'\}}{|D\Psi|^2}+C\itQ{Q^{(s)}}{g_*|\Psi|^2}.\eeq

Because $N/2\le q_0$, $2q_0'\le 2_*$, the Sobolev conjugate of $2$, so that we can use H\"older and Sobolev's inequalities ($\Psi=0$ on the boundary) to estimate the integral of $g_{*}|\Psi|^2$ over $\Og \times\{\tau\}$ by  \beqno{DPsibound} \left(\iidx{\Og\times\{\tau\}}{g_*^{q_0}}\right)^\frac1{q_0}\left(\iidx{\Og\times\{\tau\}}{|\Psi|^{2q_0'}}\right)^\frac1{q_0'}\le C\iidx{\Og\times\{\tau\}}{|D\Psi|^{2}},\quad \tau\in(0,T).\eeq
Here, we used the assumption \mref{g*bound} on $g_*$. Hence, $$\itQ{Q^{(s)}}{g_*|\Psi|^2}\le C\itQ{Q^{(s)}}{|D\Psi|^2}.$$ Using this in \mref{Psistartz} we deduce \beqno{Psistart0}\iidx{\Og\times\{s\}}{|D\Psi|^2}+\itQ{Q^{(s)}}{\llg_*|\Delta\Psi|^2}\le \iidx{\Og\times\{s'\}}{|D\Psi|^2}+C\itQ{Q^{(s)}}{|D\Psi|^2}.\eeq

Replacing $s'$ in \mref{Psistart0} by $s_k$ where $s_k$ belongs to a sequence in $(0,s)$ such that the limit of $\|D\Psi(\cdot,s_k)\|_{L^2(\Og)}$ is $\liminf_{t\to 0}\|D\Psi(\cdot,t)\|_{L^2(\Og)}$ and using \mref{limDPsi}, we then derive
\beqno{Psistart1}\iidx{\Og\times\{s\}}{|D\Psi|^2}+\itQ{Q^{(s)}}{\llg_*|\Delta\Psi|^2}\le \|D\psi\|_{L^2(\Og)}+C\int_0^s\iidx{\Og\times\{t\}}{|D\Psi|^2}\,dt.\eeq

Dropping the second term on the left side of \mref{Psistart1}, we obtain an integral  Gronwall inequality for $\|D\Psi\|_{\Og\times\{s\}}^2$ which yields $\|D\Psi\|_{\Og\times\{s\}}\le C(T)\|D\psi\|_{L^2(\Og)}$ for some constant
$C(T)$ and any $s\in(0,T)$.
This is \mref{Psibounda}. We also obtain the estimate \mref{D2Psibound} for $\Delta\Psi$ from \mref{Psistart1} and \mref{Psibounda}, using the fact that $\llg_{*}\ge \llg_0$ a positive constant.
This completes the proof of the lemma. \eproof

In addition, we have the following integrability result for $\Psi$.
\blemm{Psiboundlem} Assume as in \reflemm{Psiexistence}. Then $\Psi\in L^{\sg_N}(Q)$ with ${\sg_N}$ given by
\beqno{qPsidefa}{\sg_N}=\left\{\barr{ll} \mbox{any number in $(1,\infty)$} & \mbox{ if $N=2 $},\\ \mbox{any number in $(1,6+\frac{10}{3})$} & \mbox{ if $N= 3$},\\ \frac{2(N+2)}{N-2} & \mbox{ if $N\ge 4$}.\earr\right.\eeq
\elemm

Note that the number ${\sg_N}$ here is exactly the one defined in \mref{qPsidef} and used in \reftheo{uniweak}. The proof of this lemma based on the bounds \mref{Psibounda} and \mref{D2Psibound} and the following parabolic Sobolev imbedding inequality.

\blemm{parasobolev} Let $r^*=p/N$ if $N>p$ and $r^*$ be any number in $(0,1)$ if $N\le p$. For any sufficiently nonegative smooth functions $g,G$ and any time interval $I$ there is a constant $C$ such that\beqno{paraSobo}\itQ{\Og\times I}{g^{r^*}G^p}\le C\sup_I\left(\iidx{\Og\times\{t\}}{g}\right)^{r^*}\itQ{\Og\times I}{(|DG|^p+G^p)}\eeq   If $G=0$ on the parabolic boundary $\partial\Og\times I$ then the integral of $G^p$ over $\Og\times I$ on the right hand side can be dropped.

Furthermore, if $r<r^*$ then for any $\eg>0$ we can find a constant $C(\eg)$ such that
\beqno{paraSobo1}\itQ{\Og\times I}{g^{r}G^p}\le C\sup_I\left(\iidx{\Og\times\{t\}}{g}\right)^{r}\itQ{\Og\times I}{(\eg|DG|^p+C(\eg)G^p)}\eeq
\elemm

\bproof For any $r\in(0,1)$ and $t\in I$ we have  via H\"older's inequality
\beqno{Sobo1}\iidx{\Og}{g^rG^p}\le \left(\iidx{\Og}{g}\right)^{r}\left(\iidx{\Og}{G^\frac{p}{1-r}}\right)^{1-r}.\eeq
If $r=r^*$ then $p/(1-r)=N_*=pN/(N-p)$, the Sobolev conjugate of $p$ if $N>p$ (the case $N\le p$ is obvious),  so that the Sobolev inequality gives
$$\left(\iidx{\Og}{G^\frac{p}{1-r}}\right)^{1-r}\le \iidx{\Og}{(|DG|^p+G^p)}.$$

Using the above in \mref{Sobo1} and integrating over $I$, we  easily obtain  \mref{paraSobo}. On the other hand, if $r<r^*$, then $p/(1-r)<N_*$. A simple contradiction argument (similar to that of the proof of \reflemm{Sobointineq}) and the compactness of the imbedding of $W^{1,p}(\Og)$ into $L^{p/(1-r)}(\Og)$ imply that for any $\eg>0$ there is $C(\eg)$ such that
$$\left(\iidx{\Og}{G^\frac{p}{1-r}}\right)^{1-r}\le \eg\iidx{\Og}{|DG|^p}+C(\eg)\iidx{\Og}{G^p}.$$ We then obtain \mref{paraSobo1}. \eproof

{\bf Proof of \reflemm{Psiboundlem}:} We now have from \mref{Psibounda} and \mref{D2Psibound} that \beqno{Psibounds}  \sup_{(0,T)}\|D\Psi\|_{L^2(\Og)},\;\itQ{Q}{|\Delta\Psi|^2}<\infty.\eeq

The case $N=2$ is then obvious. Indeed, the bound for $D\Psi$ in \mref{Psibounds} and Sobolev's embedding inequality ($\Psi=0$ on $\partial\Og $) imply that $\sup_{(0,T)}\|\Psi^{q_1}\|_{L^1(\Og)}$ is finite for any $q_1\in(1,\infty)$. We need only consider $N\ge3$ below.

We first apply \reflemm{parasobolev} with $g=|D\Psi|^2$, $G=|D\Psi|$ and $p=2$.  Invoking a well known fact \cite[Corollary 9.10]{GT}, which asserts that $\|D^2 \Psi\|_{L^2(\Og)}=\|\Delta \Psi\|_{L^2(\Og)}$, we then see that
$$\itQ{Q}{|DG|^2}= \itQ{Q}{|\Delta\Psi|^2}<\infty.$$
Thus with $p=2$ and  $r_1=2/N$ as $N>2$, we get 
\beqno{DG1} \itQ{Q}{|D\Psi|^{2r_1+2}}=\itQ{Q}{g^{r_1}G^2}<\infty.\eeq

The bound for $D\Psi$ in \mref{Psibounds} and Sobolev's embedding inequality ($\Psi=0$ on $\partial\Og $) imply that $\sup_{(0,T)}\|\Psi^{q_1}\|_{L^1(\Og)}$ is finite for $q_1=\frac{2N}{N-2}$. We then apply \reflemm{parasobolev} again with $g=|\Psi|^{q_1}$, $G=|\Psi|$, and $p=2r_1+2$. Note that $p>N$ if $N=3$ only and $r^*$ can be any number in $(0,1)$ in this case. Otherwise, $r^*=p/N$.
We use \mref{DG1} to obtain an estimate for the integral of $g^rG^p=|\Psi|^{q_1r+p}$ over $Q$.  A straightforward calculation with these parameters implies that $\Psi\in L^{\sg_N}(Q)$ with ${\sg_N}=q_1r+p$ given by
$${\sg_N}=\left\{\barr{ll}  \mbox{any number in $(1,6+\frac{10}{3})$} & \mbox{ if $N= 3$},\\ \frac{2(N+2)}{N-2} & \mbox{ if $N\ge 4$}.\earr\right.$$ This is the exponent defined in \mref{qPsidefa} for $N\ge3$ and the proof is complete. \eproof

\section{Proof of Uniqueness of unbounded (very) weak solutions} \eqnoset\label{uniboundw}

In this section we present the proof of \reftheo{uniweak} on the uniqueness of very weak solutions. Let us recall its assumptions.
First of all, for some $p>2$ we assume the following continuity conditions (the definitions of $p_2,p_{\sg_N}, {\sg_N}$ are given in \mref{pr} and \mref{qPsidef}).
\bdes \item[i)]
The map $u\to \partial_u P(u)$ is continuous from $L^p(Q)$ to $L^{p_2}(Q)$,
\item[ii)] The map  $u\to \partial_u f(u)$ is continuous from $L^p(Q)$ to $L^{p_{\sg_N}}(Q)$.
\edes

We then consider very weak solutions $u$ in $L^p(Q)$ which satisfy for some $q_0\ge N/2$ \beqno{llguw} \sup_{t\in (0,T_0)}\|\hat{F}(u(t))\|_{L^{q_0}(\Og)}<\infty,\eeq
where $\hat{F}$ is a {\em convex} function in F) and satisfies $$\frac{|\partial_u f(u)|^2}{\llg(u)}\le \hat{F}(u), \quad \mbox{ for all $u\in\RR^m$}.$$

According to \refdef{wsolndef0}, we recall that a very weak solution $u$ satisfies \mref{Pwdef0} for all admissible $\fg$ with $\fg_t\in C^1(Q)$ and $\Delta \fg \in C^2(Q)$
\beqno{Pwdef}\barrl{\iidx{\Og}{\myprod{u(T),\fg(T)}-\myprod{u_0,\fg(0)}}=}{3cm}&\itQ{\Og \times(0,T)}{[\myprod{u, \fg_t}+\myprod{P(u),\Delta\fg} + \myprod{f(u),\fg}]}.\earr\eeq

{\bf Proof of \reftheo{uniweak}:} For any $u_1,u_2$ we can write  $$P(u_1)-P(u_2)=\ma(u_1,u_2)(u_1-u_2),\quad \ma(u_1,u_2):=\int_0^1 \partial_u P(su_1+(1-s)u_2)\,ds,$$ 
$$f(u_1)-f(u_2)=\mg(u_1,u_2)(u_1-u_2),\quad \mg(u_1,u_2):=\int_0^1 \partial_u f(su_1+(1-s)u_2)\,ds.$$

Using these notations, if $u_1,u_2$ are two very weak solutions with the same initial data $u_0$ then we subtract the two equations \mref{Pwdef} for $u_1,u_2$ to see that $w=u_1-u_2$ satisfies
\beqno{u1u2bded}\iidx{\Og}{\myprod{w(T),\fg(T)}}=\itQ{\Og \times(0,T)}{\myprod{w, \fg_t+\ma(u_1,u_2)^T\Delta\fg + \mg(u_1,u_2)^T\fg}}.\eeq

We consider the sequences $\{u_{1,n}\}$, $\{u_{2,n}\}$ of mollifications of $u_1, u_2$. That is, we consider $C^\infty$ functions $\eta(t)$ and $\rg (x)$ whose supports are $(-1,1)$ and $B_1(0)$ and $\|\eta\|_{L^1(\RR)} =\|\rg\|_{L^1(\RR^N)}=1$. Denote $\eta_n(t) = n\eta(t/n)$ and $\rg_n(x)=n^N\rg (x/n)$.
For $i=1,2$ define $$ u_{i,n}(t,y)=(\eta_n\fg_n)* u_i(t,y)=\int_{\RR}\iidx{\RR^N}{\eta_n (s-t)\fg_n(x-y)u_i(t,x)}\,ds.$$

For each $n$ and any given $\psi\in C^{1}(\Og)$ and $T\in(0,T_0)$ we will show in \reflemm{Psinexist} following this proof that there is a  strong solutions $\Psi_{n}$ to the systems
\beqno{Psinsys}\left\{\barr{l}\Psi_t+\ma(u_{1,n},u_{2,n})^T\Delta\Psi + \mg(u_{1,n},u_{2,n})^T\Psi=0\mbox{ in $Q:=\Og \times (0,T)$},\\\Psi=0 \mbox{ on $\partial\Og \times (0,T)$},\\\Psi(x,T) = \psi(x) \mbox{ on $\Og $}.\earr \right.\eeq
Furthermore, there is a constant $C(\|D\psi\|_{L^2(\Og)})$ such that for all $n$ \beqno{D2Wn}  \sup_{(0,T)}\|D\Psi_n\|_{L^2(\Og)},\; \|\Delta\Psi_n\|_{L^2(Q)},\;\|\Psi_n\|_{L^{\sg_N}(Q)}\le C(\|D\psi\|_{L^2(\Og)}).\eeq
One should note that the above bounds are independent of $n$ (${\sg_N}$ is defined in \mref{qPsidef}).

From the equation of $\Psi_n$ and \mref{u1u2bded} with $\fg =\Psi_n$ (this is eligible because $\Psi_n$ is a strong solution) we derive
$$\barr{lll}\iidx{\Og}{\myprod{w(T),\psi}}&=&-\itQ{Q}{\myprod{w,[\ma(u_{1,n},u_{2,n})^T-\ma(u_{1},u_{2})^T]\Delta\Psi_n}}\\&&-\itQ{Q}{\myprod{w,[\mg(u_{1,n},u_{2,n})^T-\mg(u_{1},u_{2})^T]\Psi_n}}.\earr$$ This is \beqno{keyuniw}\barr{lll}\iidx{\Og}{\myprod{w(T),\psi}}&=&-\itQ{Q}{\myprod{[\ma(u_{1,n},u_{2,n})-\ma(u_{1},u_{2})]w,\Delta\Psi_n}}\\&&-\itQ{Q}{\myprod{[\mg(u_{1,n},u_{2,n})-\mg(u_{1},u_{2})]w,\Psi_n}}.\earr\eeq

Letting $n\to\infty$, we will see that the integrals on the right hand side tend to zero. Indeed, we consider the first integral. The bound \mref{D2Wn} implies $\|\Delta\Psi_n\|_{L^2(Q)}$ is bounded uniformly in $n$ so that  we need only to show that $[\ma(u_{1,n},u_{2,n})-\ma(u_{1},u_{2})]w$ converges strongly to 0 in $L^2(Q)$. By H\"older's inequality with $q=p_2=2p/(p-2)$
$$\|[\ma(u_{1,n},u_{2,n})-\ma(u_{1},u_{2})]w\|_{L^2(Q)}\le \|\ma(u_{1,n},u_{2,n})-\ma(u_{1},u_{2})\|_{L^q(Q)}\|w\|_{L^p(Q)}.$$

As we are assuming in i) that the map $u\to \partial_u P(u)$ is continuous from $L^p(Q)$ to $L^q(Q)$ and because $u_{i,n}\to u_i$ in
$L^p(Q)$, it is clear from the definition of $\ma$ that $\ma(u_{1,n},u_{2,n})$ converges to $\ma(u_{1},u_{2})$ in $L^q(Q)$. Thus, $[\ma(u_{1,n},u_{2,n})-\ma(u_{1},u_{2})]w$ converges strongly to 0 in $L^2(Q)$. Thus, the first integral on the right hand side of \mref{keyuniw} tends to 0 as $n\to \infty$. 

Similar argument applies to the second integral to obtain the same conclusion. Using \mref{D2Wn} we see that $\|\Psi_n\|_{L^{\sg_N}(Q)}$ are uniformly bounded. We then need only to show that $[\mg(u_{1,n},u_{2,n})-\mg(u_{1},u_{2})]w$ converges strongly to 0 in $L^{{\sg_N}'}(Q)$. This is case by the assumption ii) on $\partial_u f$ and H\"older's inequality with $q=p_{\sg_N}=(\sg_N)'p/(p-(\sg_N)')$
$$\|[\mg(u_{1,n},u_{2,n})-\mg(u_{1},u_{2})]w\|_{L^{r'}(Q)}\le \|\mg(u_{1,n},u_{2,n})-\mg(u_{1},u_{2})\|_{L^q(Q)}\|w\|_{L^p(Q)}.$$

We just prove that the right hand side of \mref{keyuniw} tends to 0.
We then have $$\iidx{\Og}{\myprod{w(T),\psi}}=0 \mbox{ for any $\psi\in C^{1}(\Og)$}.$$ We conclude that $w(T)=0$ for all $T\in(0,T_0)$. Hence $u_1\equiv u_2$ on $Q$. \eproof

\brem{unboundedw} The continuity conditions i) and ii) are not needed if we discuss {\em bounded} weak solutions. Indeed, as the sequences $\{u_{1,n}\}$, $\{u_{2,n}\}$  converge to $u_1,u_2$ in $L^\infty(Q)$ wee see that $\ma(u_{1,n},u_{2,n})\to\ma(u_{1},u_{2}) $ and $\mg(u_{1,n},u_{2,n})\to\mg(u_{1},u_{2}) $ strongly in $L^\infty(Q)$.

\erem

We now provide the following lemma which establish the claim \mref{D2Wn} in the proof.
\blemm{Psinexist} For all $n$ there exist strong solutions $\Psi_n$ to the system \mref{Psinsys} and a constant $C(\|D\psi\|_{L^2(\Og)})$ independent of $n$ such that \beqno{D2Wna}  \sup_{(0,T)}\|D\Psi_n\|_{L^2(\Og)},\; \|\Delta\Psi_n\|_{L^2(Q)},\;\|\Psi_n\|_{L^{\sg_N}(Q)}\le C(\|D\psi\|_{L^2(\Og)}).\eeq
\elemm

\bproof By a change of variables $t\to T-t$, the system \mref{Psinsys} is equivalent to the following linear parabolic system for $\hat{\Psi}(x,t)=\Psi(x,T-t)$ \beqno{Psinsysa}\left\{\barr{l}\hat{\Psi}_t=\ma(u_{1,n},u_{2,n})^T\Delta\hat{\Psi} + \mg(u_{1,n},u_{2,n})^T\hat{\Psi}=0\mbox{ in $Q:=\Og \times (0,T)$},\\\hat{\Psi}(x,t)=0 \mbox{ on $\partial\Og \times (0,T)$},\\\hat{\Psi}(x,0) = \psi(x) \mbox{ on $\Og $}.\earr \right.\eeq

We will apply \reflemm{Psiexistence} here with
$\ccA = \ma^T(u_{1,n},u_{2,n})$ and $\ccG =\mg^T(u_{1,n},u_{2,n})$. We need to verify the condition AG) for such $\ccA ,\ccG$. These functions are smooth and bounded on $Q$, because $u_{1,n},u_{2,n}$ are, so that AG) is clearly satisfied with
$$\llg_* =\int_0^1 \llg(su_{1,n}+(1-s)u_{2,n})\,ds\ge \llg_0>0.$$

Thus, the existence of $\Psi_n$ is clear and we need only to establish \mref{D2Wna} by checking the condition \mref{g*bound} of \reflemm{Psiexistence}. To this end, we have to consider the function $g_*=|\ccG |^2/\llg_{*}$ and show that $\|g_*\|_{L^{q_0}(\Og)}$ is bounded for some $q_0 \ge N/2$. Here,
$$\ccG = \mg(u_{1,n}, u_{2,n})^T=\int_0^1 \partial_u f(su_{1,n}+(1-s)u_{2,n})^T\,ds$$

Denote $w(s)=su_{1,n}+(1-s)u_{2,n}$. Writing $|\partial_u f(w(s))|=\frac{|\partial_u f(w(s))|}{\llg^\frac12(w(s))}\llg^\frac12(w(s))$ and using H\"older's inequality, we have
$$\int_0^1|\partial_u f(w(s))|ds\le \left(\int_0^1\frac{|\partial_u f(w(s))|^2}{\llg(w(s))}ds\right)^\frac12\left(\int_0^1\llg(w(s))ds\right)^\frac12.$$
This implies for each $t\in(0,T)$ that $$\frac{|\ccG |^2}{\llg_*}\le \int_0^1\frac{|\partial_u f(w(s))|^2}{\llg(w(s))}ds \le \int_0^1 \hat{F}(s)\,ds,$$
where $\hat{F}$ is the convex function specified in F) of the theorem and satisfying $\frac{|\partial_u f(w)|^2}{\llg(w)}\le \hat{F}(w)$. Therefore,
$$\left\|\frac{|\ccG |^2}{\llg_*}\right\|_{L^{q_0}(\Og)}\le \int_0^1 \left\|\hat{F}(su_{1,n}+(1-s)u_{2,n})\right\|_{L^{q_0}(\Og)}ds.$$

As $\hat{F}$ is convex, by Jensen's inequality $\hat{F}(su_{1,n}+(1-s)u_{2,n})\le s\hat{F} (u_{1,n})+(1-s)\hat{F} (u_{2,n})$. Similarly, because $\|\eta_n\rg_n\|_{L^1(\RR^N)}=1$, for $i=1,2$ $\hat{F} (u_{i,n})\le (\eta_n\rg_n) *\hat{F} (u_i)$ so that
$$\|\hat{F} (u_{i,n}(t))\|_{L^{q_0}(\Og)} \le \int_\RR \eta_n(s-t)\|\rg_n *_x \hat{F}(u_{i}(t))\|_{L^{q_0}(\Og)}\,ds.$$
Here, $*_x$ denotes the convolution in $\RR^N$. Obviously, $\|\rg_n *_x \hat{F}(u_{i}(t))\|_{L^{q_0}(\Og)}\le\|\hat{F} (u_{i}(t))\|_{L^{q_0}(\Og)}$. As $\|\eta_n\|_{L^1(\RR)}=1$, we find a constant $C_0$, by the assumption \mref{llguw}, such that
\beqno{llgun} \|\hat{F} (u_{i,n}(t))\|_{L^{q_0}(\Og)} \le \|\hat{F} (u_{i}(t))\|_{L^{q_0}(\Og)}\le C_0 \mbox{ for all $t\in(0,T)$ and integer $n$}.\eeq

Hence,  $\||\ccG|^2\llg_{*}^{-1} \|_{L^{q_0}(\Og)}\le C_0$   on $(0,T)$ for some $q_0\ge N/2$ and we can apply  \reflemm{Psinexist} to obtain a constant $C(\|D\psi\|_{L^2(\Og)})$, which is also dependent of $C_0$ and $T$ but independent of $n$ such that $$\sup_{(0,T)}\|D\hat{\Psi}_n\|_{L^2(\Og)},\; \|\Delta\hat{\Psi}_n\|_{L^2(Q)}\le C(\|D\psi\|_{L^2(\Og)}).$$
Because $\Psi_n(x,t)=\hat{\Psi}_n(x,T-t)$, the above estimate implies the bound for the first two terms in \mref{D2Wna}. The bound for the last term in \mref{D2Wna} comes from  \reflemm{Psiboundlem}.  \eproof

In many models in application, the components of the maps $P,f$ of \mref{ep1} are polynomials of $u\in\RR^m$. Our uniqueness theorem then applies to these models to show that unbounded weak solutions are unique if they satisfy sufficient integrability as in the following

\bcoro{ALqcont} For some $k,l>0$ with $2l-k\ge1$ assume that $\partial_u P(u)$ and $\partial_u f(u)$ have polynomial growths
$$|\partial_u P(u)|\le C(|u|^k+1),\; |\partial_u f(u)|\le C(|u|^l+1).$$ Then the uniqueness conclusion of \reftheo{uniweak} applies to weak solutions in the space $L^p(Q)\cap L^\infty((0,T_0),L^r(\Og))$ if $p\ge\max\{2(1+k),(\sg_N)'(1+l)\}$ and $r\ge (2l-k)N/2$.
\ecoro

\bproof We need to verify first the the assumptions i) and ii) of \reftheo{uniweak}. It is clear from its proof that we need only to establish the convergences $\ma(u_{1,n},u_{2,n})\to \ma(u_{1},u_{2})$ and $\mg(u_{1,n},u_{2,n})\to \mg(u_{1},u_{2})$ in $L^q(Q)$ for appropriate $q$'s along some subsequences of $\{u_{1,n}\}$,$\{u_{2,n}\}$ which converge to $u_1$,$u_2$ in $L^p(Q)$. In particular,  by the Riesz-Fisher theorem we can find subsequences of $\{u_{i,n}\}$ and functions $\hat{u}_i\in L^p(Q)$ such that, after relabeling, $u_{i,n}\to u_i$ and $|u_{i,n}|\le \hat{u}_i$  a.e. in $Q$. The growth condition of $\partial_u P(u)$  then implies $|\ma(u_{1,n},u_{2,n})|\le |\hat{u}_1|^k+|\hat{u}_2|^k+1$, a function in $L^{p/k}(Q)$. Furthermore, $\ma(u_{1,n},u_{2,n})\to \ma(u_{1},u_{2})$ a.e. in $Q$ because $\ma$ is continuous. By the Dominated convergence theorem, we see that $\ma(u_{1,n},u_{2,n})\to \ma(u_{1},u_{2})$ in $L^{p/k}(Q)$. This also yields the convergence in $L^q(Q)$ for $q=p_2=2p/(p-2)$ because $q\le p/k$ (as $p\ge 2(1+k)$). Similarly, from the growth condition of $\partial_u f$, we replace $2,k$ respectively by $(\sg_N)',l$ in the argument to see that
$\mg(u_{1,n},u_{2,n})\to \mg(u_{1},u_{2})$ in $L^{p/l}(Q)$ and thus in $L^q(Q)$ for $q=p_{\sg_N}$.

Finally, from the growth assumptions on $\partial_u P$ and $\partial_u f$, we see that $|\partial_u f(u)|^2/\llg (u)\le C|u|^{2l-k}$ for $|u|$ large so that we can take the function $\hat{F}$ defined in F) to be $C|u|^{2l-k}$. Because $2l-k\ge1$, $\hat{F}$ is convex. We then have  $|\hat{F}(u)|^{q_0}\le C|u|^{(2l-k)q_0}$. Because $u\in L^\infty((0,T_0),L^r(\Og))$ for some $r \ge (2l-k)N/2$ from the assumption of the corollary,  we have $\sup_{(0,T_0)}\|\hat{F} (u)\|_{L^{q_0}(\Og)}$ is finite for some $q_0\ge N/2$.  We see that all assumptions of of \reftheo{uniweak} are verified here. This completes the proof.
\eproof

\brem{Pubounded} If $P_u(u)$ is bounded, i.e., $k=0$ then we can allow $p=2$ in \reftheo{uniweak} (the condition i) is then dropped) and the above corollary. Indeed, from the proof of the theorem, we need the sequence $\{[\ma(u_{1,n},u_{2,n})-\ma(u_{1},u_{2})]w\}$, $w=u_1-u_2$, converges strongly in $L^2(Q)$. But this is obvious because this sequence converges pointwise in $Q$ and is bounded by $C|w|$, a function in $L^2(Q)$. The Dominated convergence theorem applies again to give the corollary.
\erem

As we mentioned in \refsec{mainres}, our \refdef{wsolndef} is the most general one as one needs at least that the derivatives of the test function $\fg $ are defined and bounded so that we need only $\fg \in C^1(Q)$ so that we need only that $u\in L^\infty((0,T_0),L^1(\Og))$ and $D(P(u))\in L^1(Q)$. This assumption is too weak in order to verify the integrability conditions of \reftheo{uniweak}.

On the other hand, if we allow more general test function $\fg$ then the space of weak solutions will be smaller and the uniqueness result can be applied easily and almost immediate in some cases. As an example, we will discuss an application of \refcoro{ALqcont} to the SKT system and its generalizations.
We now  prove \refcoro{SKTuni} concerning the generalized solutions from $V^2(Q)$.

{\bf Proof of \refcoro{SKTuni}:} Let $u$ be a generalized solution from $V^2(Q)$.
Now, the space of admissible test functions is $W_2^{1,1}(Q)$ and it is clear that in order for the integrals in \mref{wdef} to be finite for all $\fg \in W_2^{1,1}(Q)$ we should assume further that $u\in L^\infty((0,T),L^2(\Og))$ and $D(P(u))\in L^2(Q)$. The first assumption is satisfied because  $u\in V^2(Q)$. We now let $g=|u|$, $G=|P(u)|$ and $p=2$ in \mref{paraSobo}. As $|P(u)|\sim |u|^{k+1}$, we see that $u\in L^{2r^*+2k+2}(Q)$ for some $r^*>0$. The condition in \refcoro{ALqcont} that $u\in L^{p}(Q)$ with $p\ge 2+2k$ is then obvious and this also implies $p\ge (\sg_N)'(1+l)$ as $l=k$ and $(\sg_N)'<2$. 
From \refcoro{ALqcont}, we need only that  $\sup_{(0,T_0)}\|u\|_{L^r(\Og)}$ is finite for some $r\ge kN/2$. This condition is clearly satisfied for generalized solutions ($r=2$) to the usual SKT system, where $P(u)$ has quadratic growth in $u$ (so that $k=1$), in domains with dimension $N\le4$. Thus, $u$ is unique. \eproof

The above proof also implies
\bcoro{genSKTuni} Generalized solutions to generalized SKT systems are unique as long as $1\le k\le 4/N$. \ecoro

\section{Regularity} \eqnoset \label{regsec}

In this section, we consider the regularity problem of weak solutions of \mref{ep1} when $\Og$ is a planar domain, i.e., $N=2$. The main idea of the proof is simple. We will show that there exists a strong solution $u$ to \mref{ep1}. This solution is of course a weak one. By the uniqueness result for weak solutions in the previous section, any weak solution satisfying sufficient integrability is then exactly this strong solution $u$ and is in fact classical.

We first establish the existence of strong solutions of \mref{ep1}. To this end, we embed \mref{ep1} in the following family of systems with $\sg\in[0,1]$
\beqno{mainparaBfam}\left\{\barr{l} w_t-\Delta(P(w))= \sg^2 f(w),\mbox{ $(x,t)\in \Og\times(0,T_0)$},\\\mbox{$w=0$ on $\partial \Og\times(0,T_0)$},\\ w(x,0)=\sg u_{0}(x) ,\quad x\in \Og. \earr\right.\eeq
We then assume the following main apriori integrability condition on \mref{mainparaBfam}.
\bdes \item[S)] There are $q_0>1$ and a constant $C_1$, which may depends on $T_0$ but independent of $\sg \in[0,1]$, such that any {\em strong} solutions $w$ of this \mref{mainparaBfam} satisfy 
\beqno{famass} \sup_{(0,T_0)}\|\llg (w)\|_{L^{q_0}(\Og\times\{t\})},\; \sup_{(0,T_0)}\|w\|_{L^{q_0}(\Og\times\{t\})}\le C_1.\eeq

Furthermore, $\llg(u), f(u)$ have polynomial growths in $|u|$ and for some constant $C$ all all $u\in \RR^m$ \beqno{fllg} |\llg_u(u)||u|\le C\llg(u) \mbox{ and } |f(u)|\le C(1+|u|)(1+\llg(u)).\eeq

\edes

We should note that the inegrability condition \mref{famass} is a very mild one, especially it needs only be satisfied for {\em strong solutions}, and can be verified in many models.

We first establish the existence of a strong solution to \mref{ep1}. To this end, we will use the theory in \cite{dlebook} which needs smoother initial data $u_0$.
\bprop{strongthm} Assume S) and $u_0\in W^{1,p}_0(\Og)$ for some $p>2$.
Then there is a strong (classical) solution to \mref{ep1}.

\eprop

\bproof The system \mref{ep1} can be written as, with $A(u)=P_u(u)$ \beqno{mainparaBfamz}\left\{\barr{l}u_t-\Div(A(u)Du)=f(u),\quad (x,t)\in \Og\times(0,T_0),\\ u=0, \quad (x,t)\in \partial\Og\times(0,T_0),\\ u(x,0)=u_0(x)\quad x\in\Og.\earr\right.\eeq

To establish the existence of a strong solution to this system, we apply \cite[Theorem 3.4.1]{dlebook} here by verifying its assumptions. First of all, we need to show that the number $\mathbf{\LLg}=\sup_{u\in\RR^m}\LLg(u)$,  with $\LLg(u)=|\llg_u(u)|/\llg(u)$, is finite.  Since $\llg(u)\ge \llg_{0}>0$, if $|u|$ is bounded then so is $\LLg(u)$. For large $|u|$ we use the assumption in \mref{famass} that $|\llg_u(u)|\le C\llg(u)/|u|$ to see that $\LLg(u)\le C/|u|$ is also bounded. Hence, the number $\mathbf{\LLg}$ is finite.

Next, also following \cite[Theorem 3.4.1]{dlebook}, we  embed \mref{mainparaBfamz}  in the following family 
\beqno{gensysfamB}\left\{\barr{l}u_t-\Div(A(\sg u)Du)=\sg f(\sg u),\quad (x,t)\in \Og\times(0,T_0),\\ u=0, \quad (x,t)\in \partial\Og\times(0,T_0),\\ u(x,0)=\sg u_0(x)\quad x\in\Og.\earr\right.\eeq 

The most non trivial condition of the theory in \cite{dlebook} needs to be checked is that the strong solutions of this family has small BMO norm in small ball uniformly. Namely,

\bdes\item[(Sbmo)] For any given $\mu>0$ there is $R>0$ depending only on $\mu$ and the parameters of the system \mref{ep1} (but {\em not} on $\sg$) such that any strong solution $u$ to the family \mref{gensysfamB} satisfies:
for any ball $B_R$ in $\RR^N$ with $\Og_R=B_R\cap \Og\ne\emptyset$ $$\sup_{t\in(0,T_0)}\|u(\cdot,t)\|_{BMO(\Og_R)} \le \mu.$$
\edes

We is going to verify this property. Multiplying $\sg>0$ to the equation in \mref{gensysfamB}, we see easily that $w=\sg u$ is a strong solution to \mref{mainparaBfam}.

For $\sg=1$ and a strong solution $w$ of \mref{mainparaBfam} we can multiply the system of $w$ by $P(w)_t$ and follow the proof of \cite[Lemma 5.3.2]{dlebook} and use the fact that $$\frac{d}{dt}|D(P(u))|^2=2\myprod{DP(u),(D(P(u))_t}=2\myprod{DP(u),D(P(u)_t)},$$ to prove that: There is an absolute constant $C^*$ such that
or any $t\in(0,T_0)$ \beqno{Aeqn}\iidx{\Og\times\{t\}}{\llg(w)|w_t|^2} +\frac{d }{d t}\iidx{\Og\times\{t\}}{|A(w)Dw|^2}\le C^*\iidx{\Og\times\{t\}}{\llg(w)|f(w)|^2}.\eeq 

This yields a Gronwall inequality for $\|A(w)Dw\|_{L^2(\Og)}^2$. Indeed, as in the proof of \cite[Proposition 5.3.1]{dlebook} (see also \reftheo{wsequal1} in this paper) for general $N$,  under the growth condition \mref{fllg} in S) and its assumption that there is a constant $C_1$ such that
\beqno{keyllgf}\sup_{(0,T_0)}\|\llg (w)\|_{L^{q_0}(\Og)},\; \sup_{(0,T_0)}\|w\|_{L^{q_1}(\Og)}\le C_1\eeq for some $q_0>N/2$ and $q_1>2N/(N+2)$, which is fulfilled here by \mref{famass} as $N=2$, we can prove that there is a constant $C$ such that\beqno{llgf} \iidx{\Og\times\{t\}}{\llg(w)|f(w)|^2}\le C\iidx{\Og\times\{t\}}{|A(w)Dw|^2}+C(C_1).\eeq This and \mref{Aeqn} imply
\beqno{Aeqn0}\frac{d }{d t}\iidx{\Og\times\{t\}}{|A(w)Dw|^2}\le C^*\left[C\iidx{\Og\times\{t\}}{|A(w)Dw|^2}+C(C_1)\right].\eeq

On the other hand, we can apply \reflemm{limDuat0} (and \refrem{LpDADBrem0}) to the system \mref{mainparaBfamz} (in divergence form) here, with
$$\ccA = A(w),\; \ccG = \int_0^1 \partial_u f(sw)ds$$
to see that $\liminf_{t\to0}\|Dw\|_{L^2(\Og)} \le \|Du_0\|_{L^2(\Og)}$. The condition AG) is satisfied here because $w$ is a strong solution and continuous at $t=0$. Also, by \refrem{LpDADBrem0}, we do not need the boundedness of $D\ccA = DA(w)$ here for systems in divergence form like \mref{gensysfamB}. It then follows 
$$\liminf_{t\to0}\|A(w)Dw\|_{L^2(\Og)} \le \| A(u_0)\|_{L^\infty(\Og)}\|Du_0\|_{L^2(\Og)}.$$ Therefore, from the Gronwall inequality \mref{Aeqn0} for $\|A(w)Dw\|_{L(\Og)}^2$ and the fact that $\llg(w)|Dw|$ is comparable to $|A(w)Dw|$, we have
$$\sup_{t\in[0,T_0]}\iidx{\Og}{\llg^2(w)|Dw|^2} \le \| A(u_0)\|_{L^\infty(\Og)}^2\|Du_0\|_{L^2(\Og)}^2+C^*C(C_1),$$ where $C_1$ is the constant in \mref{famass}.

For $\sg\ne1$ we replace $f(w),u_0$ respectively by $\sg^2 f(w),\sg u_0$ in the above argument. The constant $C^*$ in \mref{Aeqn} and the above estimates is now $\sg^2 C^*$ accordingly. As $w=\sg u$, the above argument shows that
$$\sup_{t\in[0,T_0]}\iidx{\Og}{\llg(\sg u)^{2}|D(\sg u)|^2} \le \sg^2C(\|Du_0\|_{L^2(\Og)},C_1).$$

Using the fact that $\llg (\sg u)\ge\llg_0$, we obtain $$\llg_{0}^{2}\sup_{t\in[0,T_0]}\iidx{\Og}{|D\sg u|^2} \le \sg^2C(\|Du_0\|_{L^2(\Og)},C_1).$$ Thus, the strong solutions $u$ of \mref{gensysfamB} satisfy \beqno{du2Omest}\sup_{t\in[0,T_0]}\iidx{\Og}{|Du|^2} \le \llg_{0}^{-2}C(\|Du_0\|_{L^2(\Og)},C_1).\eeq

As $N=2$, a simple use of Poincar\'e's inequality, the continuity of integral  and the above estimate show that strong solutions $u$ to \mref{gensysfamB} satisfy the (Sbmo) condition in \cite{dlebook} uniformly in $\sg\in(0,1]$ (see also \cite[Corollary 3.4.4]{dlebook} for more details on the implication of the property Sbmo) from \mref{du2Omest}). 

By Sobolev's inequality and because $u=0$ on the boundary, \mref{du2Omest}  implies that $\|u\|_{L^q(\Og)}$ is uniformly bounded for any $q\ge1$. 
From the polynomial growths of $\llg$ and $f$, we now see that $\llg (u)$, $|f(u)|\llg^{-1}(u)$ are in bounded by powers of $|u|$ so that their integrability conditions in \cite[Theorem 3.4.1]{dlebook} are verified. The last condition needs to be checked is 
\beqno{llgDu2bound1} \int_{0}^{T_0}\iidx{\Og}{|Du|^2}dt\le C(T_0)\eeq
for some constant $C(T_0)$. But this is an immediate consequence of \mref{du2Omest}.

We thus verified all conditions of \cite[Theorem 3.4.1]{dlebook} and therefore obtain the existence of a strong solution of \mref{ep1}. \eproof

We now turn to the regularity problem of unbounded weak solutions. 
The following theorem shows that unbounded weak solutions are in fact smooth if they satisfy sufficient higher integrability. We would like to emphasize that the common and crucial assumption on the smallness of their BMO norms in small balls is not needed here. 

\btheo{wsequal} Assume the conditions of \reftheo{uniweak} and S). If $u$ is a very weak solution of \mref{ep1}, with initial data $u_0\in W^{1,p_0}(\Og)$ for some $p_0>2$, in the sense of \refdef{wsolndef0},  and satisfies the integrability conditions of \reftheo{uniweak} then $u$ is a classical one. \etheo
\bproof By S), we already show in \refprop{strongthm} that there exists a strong solution $u$ to \mref{ep1}. This solution is of course a {\em bounded} weak one and satisfies the integrability of \reftheo{uniweak}. Now, any (very) weak solution satisfying sufficient integrability \reftheo{uniweak} must be this strong solution $u$ by the uniqueness. Thus, these weak solutions are in fact classical and the theorem is proved. \eproof

\brem{regrem} Even for systems like \mref{ep1} on planar domains, our regularity result asserted in \reftheo{wsequal} is remarkable because the direct approach in literature (e.g., see \cite{GiaS}) cannot apply here as the basic Caccioppoli and Poincar\'e inequalities are not available for weak solutions in the sense of \refdef{wsolndef0} (one cannot test the systems \mref{Pwdef} with $u$, which is not admissible, to obtain such inequalities). In addition, the key estimate \mref{du2Omest} is not available for weak solutions. In this proposition, we worked with {\em strong} solutions so that the integrals in starting inequality \mref{Aeqn} (and those follow) were all finite in order to derive \mref{du2Omest}. \erem 

Let us consider an alternate version of S) and connect the condition sets of \reftheo{uniweak} and \refprop{strongthm}. Following \reftheo{uniweak} and with a slight abuse of notations, we define (without the convexity assumption in F))
\beqno{hatFdef}\hat{F}(u):= |\partial_u f(u)|^2\llg^{-1}(u).\eeq

\bdes \item[S')] There are $q_0\ge1$, $\bg\in (0,1)$ and constants $C_0,C_1$, which may depends on $T_0$ but independent of $\sg \in[0,1]$, such that any {\em strong} solutions $w$ of \mref{mainparaBfam} satisfies 
\beqno{famass1}  \sup_{(0,T_0)}\|\hat{F}(w)\|_{L^{q_0}(\Og\times\{t\})}\le C_0,\; \sup_{(0,T_0)}\|(\llg (w)|w|)^\bg\|_{L^{1}(\Og\times\{t\})}\le C_1.\eeq

Furthermore, $\llg(u), f(u)$ have polynomial growths in $|u|$ and for some positive constants $C,M$ if $|u|\ge M$ then  \beqno{fllg1} |\llg_u(u)||u|\le C\llg(u), \eeq \beqno{fllg1a} |f(u)|\le C|u||\partial_u f(u)|.\eeq

\edes

\btheo{wsequal1} The conclusion of \reftheo{wsequal} still holds if S) is replaced by S'). \etheo

\bproof We just need to show that the condition S') will lead to the same Gronwall inequality \mref{Aeqn0}  for $\|A(w)Dw\|_{L(\Og\times\{t\})}^2$ so that the proof of \refprop{strongthm} can continue to provide the existence of a unique strong solution. Together with the uniqueness of \reftheo{wsequal} we obtain our assertion here.

We can argue as in the proof of \refprop{strongthm} until we get \mref{llgf}.  We need to estimate the last integral of $\llg(w)|f(w)|^2$ in this inequality to obtain \mref{Aeqn0}.

To this end,  for each for $t\in(0,T)$ we denote $\Og_{M,t}:=\{(x,t)\,:\, |w(x,t)|\ge M\}$. On $\Og\times\{t\}\setminus \Og_{M,t}$, $|w|<M$ so that the integrand $\llg(w)|f(w)|^2$ is bounded by some constant depending on $M$. On $\Og_{M,t}$,  we use the condition  \mref{fllg1a} and the notation \mref{hatFdef} to see that $\llg(w)|f(w)|^2\le C(\llg(w)|w|)^2|\partial_u f(w)|^2\llg^{-1}(w)\le C(\llg(w)|w|)^2\hat{F}(w)$. Thus, we will estimate the integral of $(\llg(w)|w|)^2\hat{F}(w)$ below. If $q_0\ge N/2$ then we can find $q\in [1,2_*]$ such that $N/2 \le (q/2)'\le q_0$ and apply the H\"older inequality and the \mref{famass1} to have
$$\barr{lll}\iidx{\Og_{M,t}}{(\llg(w)|w|)^2\hat{F}(w)}&\le& \left(\iidx{\Og_{M,t}}{(\llg(w)|w|)^q}\right)^\frac{2}{q}\|\hat{F}(w)\|_{L^{(\frac{q}{2})'}(\Og)}\\&\le& C_0\left(\iidx{\Og_{M,t}}{(\llg(w)|w|)^q}\right)^\frac{2}{q}.\earr$$

Because $q\le 2_*$, we can apply \reflemm{Sobointineq} with $p=2$ and $W=\llg(w)|w|$ to have
$$\iidx{\Og_{M,t}}{(\llg(w)|w|)^q} \le \iidx{\Og_{M,t}}{|D(\llg(w)|w|)|^2} + C(\bg)\left(\iidx{\Og_{M,t}}{(\llg(w)|w|)^\bg}\right)^\frac{2}{\bg}$$ where $\bg\in (0,1)$. Let $\bg$ be the exponent in \mref{famass1}. The last integral is then bounded by $C_1$.  On the other hand, by \mref{fllg1} we have $|\llg_u(w)||w|\le C\llg(w)$ so that$$ |D(\llg(w)|w|)|\le \llg(w)|Dw|+ |\llg_u(w)||w||Dw|\le C\llg(w)|Dw|.$$

Putting these estimates together and using $|\llg(w)Dw|\sim|A(w)Dw|$, we see that
\beqno{llgf0} \iidx{\Og\times\{t\}}{\llg(w)|f(w)|^2}\le C\iidx{\Og\times\{t\}}{|A(w)Dw|^2}+C(C_0,C_1)+C(M).\eeq

Hence, the Gronwall inequality \mref{Aeqn0}  for $\|A(w)Dw\|_{L(\Og)}^2$ continues to hold and the proof of the theorem can go on as before. \eproof

\brem{llgfgrowth} The assumption \mref{fllg1} implicitly implies that $\llg(u)$ must have a polynomial growth in $|u|$. Meanwhile, \mref{fllg1a} does not require such growth for $f(u)$. However, as we see in the proof of \refprop{strongthm}, we need that $\|f(u)\llg^{-1}(u)\|_{L^q(\Og)}$ is bounded from some large $q$ so that a polynomial growth for $f(u)$ seems to be necessary.  
\erem

\brem{Nrem} Our argument shows that the  Gronwall inequality \mref{Aeqn0} holds under S') with $q_0\ge N/2$ but eventually we still have to assume $N=2$ in order for the property Sbmo) can be verified in accordance with the theory in \cite{dlebook} to establish the existence of strong solutions. \erem

We conclude the paper by presenting some examples. In particular, we consider \mref{ep1} with polynomial growth for its data as in \refcoro{ALqcont} and for simplicity we assume that $k=l$. that is $P,f$ have the same growth $k+1$. Immediately, we have the following result.

\bcoro{genSKTreg} Assume the growth conditions as in \refcoro{ALqcont} with $k=l$. Assume also S) (or S')). If $u$ is a weak solution of \mref{ep1} in $L^p(Q)\cap L^\infty((0,T_0),L^r(\Og))$ for some $p\ge 2(1+k)$ and $r\ge k$ then $u$ is classical. \ecoro
 
\bproof The proof is almost obvious. We apply \reftheo{wsequal} (or \reftheo{wsequal}) with \refcoro{ALqcont} in place of \reftheo{uniweak}. The integrability condition of $u$ in \refcoro{ALqcont} is already assumed here, noting that ${\sg_N}>2$ so that $(\sg_N)'<2$ and the condition $p\ge 2(1+k)$ alone is sufficient. \eproof

The stated \refcoro{SKTrek} in \refsec{mainres} on the regularity of generalized weak solutions to the SKT systems now follows easily. 

{\bf Proof of \refcoro{SKTrek}:}  
We recall the assumption \mref{SKTfu} \beqno{SKTfu0} \myprod{f(w),w}\le \eg_0\llg(w)|w|^2+ C|w|^2 \eeq for some positive constants $C,\eg_0$. We will show that if $\eg_0$ is sufficiently small \refcoro{genSKTreg} applies to give that the generalized solutions from $V_2(Q)$ of  of the SKT system and its generalized versions (with $k<2$) on planar domains are classical.

The description of the SKT system and its generalized versions clearly implies the growth condition \mref{fllg} of S). We will show that  the integrability condition \mref{famass} also holds under \mref{SKTfu} if $\eg_0$ is sufficiently small.

Let $w$ be a strong solution of \mref{mainparaBfam}. Testing the system with $w$, we obtain for any $T\in (0,T_0)$ and $Q_T=\Og\times(0,T)$
$$ \sup_{t\in(0,T)}\iidx{\Og}{|w|^2} +\itQ{Q_T}{\llg(w)|Dw|^2}\le C\sg^2\itQ{Q_T}{\myprod{f(w),w}}+\iidx{\Og}{|u_0|^2}.$$

From the assumption \mref{SKTfu} and the polynomial growth  $\llg(w)\sim (\llg_0+|w|)^k$ we deduce
$$ \sup_{t\in(0,T)}\iidx{\Og}{|w|^2} +\itQ{Q_T}{|w|^k|Dw|^2}\le C\sg^2\itQ{Q_T}{\eg_0|w|^{k+2}+|w|^{2}} +C.$$

Because $w=0$ on $\partial\Og$ and $N=2$, using the Poincar\'e inequality, we can find a constant $C$ depending on the diameter of $\Og$ such that $$\iidx{\Og}{|w|^{k+2}}\le C\iidx{\Og}{|w|^{k}|Dw|^2}.$$ Hence, if $\eg_0$ is sufficiently small then we can deduce from the above two inequalities that
$$ \sup_{t\in(0,T)}\iidx{\Og}{|w|^2} \le C\int_0^T\iidx{\Og}{|w|^2}dt +C.$$
This is a Gronwall inequality for $\|w\|_{L^2(\Og)}^2$ and yields a bound for $\|w\|_{L^2(\Og)}$. As we are assuming $\llg(u)\sim (\llg_0+|u|)^k$ and $k<2$, the integrability condition \mref{famass} is verified for some $q_0>1$. The assertion then follows as in \refcoro{genSKTreg}. and the proof is complete. \eproof

\bibliographystyle{plain}

\end{document}